# A modified greedy algorithm to improve bounds for the vertex cover number

*R. Dharmarajan[1, 2] and D. Ramachandran[1]

[1]Niels Abel Foundation, Palakkad 678011, Kerala, India.
claudebergedr@gmail.com
[2]Department of Mathematics, Amrita School of Engineering - Coimbatore,
Amrita Vishwa Vidyapeetham, India.
r_dharmarajan@cb.amrita.edu

*Corresponding author

## Abstract

In any attempt at designing an efficient algorithm for the minimum vertex cover problem, obtaining good upper and lower bounds for the vertex cover number could be crucial. In this article we present a modified greedy algorithm of worst-case time complexity $O(n^3)$ to obtain bounds for the vertex cover number of an input graph of order n. Using simple facts, the proposed algorithm computes a lower bound for the vertex cover number. Then using this lower bound it outputs a minimal vertex cover and hence gives an upper bound. The algorithm ensures the output vertex cover is always minimal, which feature is an improvement upon the existing greedy algorithms.

**Keywords:** Degree; neighbour; vertex cover; independent set; lower bound; upper bound.

**2010 Mathematics Subject Classification:** 05C69, 05C70

## 1. Introduction

The minimum vertex cover problem (MVC) is to find a vertex cover of the smallest possible size in a given graph. It is an NP-complete problem [7, 16] and is computationally equivalent to the Maximum Independent Set problem. There are exact algorithms for the MVC [1, 9] but all these take exponential time and so are not suited to practical use in large graphs. Some improvements in (exponential) run time have been achieved [5, 6] but an efficient exact algorithm for the MVC is lacking. Since the MVC is NP-hard [11], no polynomial time algorithm to solve the problem is expected to be found.

Nevertheless, it is worthwhile to attempt algorithms for a few good reasons. The first is that the MVC has important applications in many domains [13, 14] - for instance, VLSI design, Bio-informatics and networking. Its real-world applications include the de novo genome assembly, finding phylogenetic trees based on protein domain information, resolving conflicts in gene alignments, reconfiguring memory arrays, placing ATM's in a city and keyword-based text summarization process.

The second reason is that new algorithms that improve some aspects of existing ones are always desirable since they may handle practical instances with increased efficiency and /or improved proximity to optimality. Third, modified algorithms could show new directions. Indeed, quite a few algorithms are modifications of earlier ones. For instance, Greedy, Maximum Degree Greedy [18] and Greedy Independent Cover [2] are very popular algorithms for the MVC and are modifications of the greedy algorithm. Edge Deletion, Depth First Search, List Left (LL), Sorted LL and Anti Sorted LL are other popular modified algorithms [2] for the MVC.

Modified algorithms do not guarantee optimal solutions but in some instances return solutions that are `close enough' to optimal. Some of them do return optimal solutions on some graph classes but fail to replicate the success elsewhere. For instance, the Greedy Independent Cover always returns an optimal cover on any tree [2] but does not perform as well with other graph classes.

Trading optimality for efficiency brought some more success for some graph classes [2, 19, 21]. Such success included algorithms with improved approximation ratios. So far the best known approximation ratio for the MVC is achieved by Karakostas' solution [15].

In general, exact algorithms return optimal solutions but cannot run in polynomial time on every graph while on the other hand there are modified algorithms that run in polynomial time but fall short of optimal solutions in overwhelmingly many cases. Still, the latter algorithms continue to attract researchers' attention because they can handle specific instances rather efficiently, especially in applications where optimal solutions are not a must and solutions `close enough to optimal' suffice [21].

The MVC has also generated considerable research for algorithm-friendly criteria simpler than existing ones. This has led to searches for closer bounds for the vertex cover number [9]. This article is one such attempt, and the central problem here is efficient computing of lower and upper bounds for the vertex cover number of a given graph. The main motivations for the algorithm presented in this article lie in the Maximum Degree Greedy.

The rest of this article is organized in sections 2 through 7. Section 2 recalls some definitions and standard results of graph theory. Section 3 is devoted to the theoretical facts on which the proposed algorithm (named Enter-Exit Greedy Algorithm) is built. Section 4 outlines the algorithm and follows up by establishing relevant conclusions. Section 5 is devoted to the time complexity of the algorithm. In section 6 there are discussions on the algorithm, including how it compares with earlier algorithms. This section also contains an illustration of the running of the algorithm on a graph. Concluding remarks form section 7.

## 2. Basic theory

Most of the definitions and notation given in this section are from [3, 4, 10, 12, 20]. Here they are meant for ready reference.

**2.1: Graph**

If V is a set, then $2^V$ denotes the *power set* of V – that is, the set of all the subsets of V (including the empty set $\phi$); and $2^{V*}$ denotes the set of all nonempty subsets of V – that is, $2^{V*} = 2^V -\{\phi\}$. The *cardinality* (or, *size*) of a finite set V is denoted by $|V|$, and is the number of elements in V. \\

A *simple undirected graph* is an ordered pair G = (V, E) $ where V is a nonempty finite set and E $\subset 2^{V*}$ such that (i) $|X| \leq 2$ for each X ∈ E and (ii) $\cup_{X \in E} X \subset V$. The sets V and E are, respectively, the *vertex set* and the *edge set* of the graph G. Each element of V is a *vertex* of G and each member of E is an *edge* of $G$. The integers $|V|$ and $|E|$ are, respectively, the *order* (= the number of vertices) and the number of edges of G. A *loop* is an edge X with $|X| = 1$. G is *loop-free* if $|X| = 2$ for each X ∈ E. If G is loop-free and { x, y} is an edge of G then x and y are the *end points* (or, *ends*) of this edge.

All the graphs in this article are undirected, simple and loop-free. Let G = (V, E) be a graph. The expressions x ∈ V and x ∈ G will both mean x is a vertex of G; similarly, {x, y} ∈ E and {x, y} ∈ G will both mean {x, y} is an edge of G. In subsections 2.2 through 2.5, G = (V, E) is assumed.

### 2.2: Adjacency, degree and neighbourhood

Two distinct vertices x and y are *adjacent* if {x, y} ∈ E; then x and y are *neighbours* in G. The set N(x) = {y ∈ V : x and y are adjacent} is the *neighbourhood* (or, *open neighbourhood*) of the vertex x. The *degree* of a vertex x in G is denoted by dx (or, dx(G)) and is defined as dx = $|N(x)|$. If A is a nonempty subset of V and x ∈ V, then the *degree of x restricted to* A is denoted by dx(A) and is defined as dx(A) = $|N(x) \cap A|$.

### 2.3: Connectedness

A *path* in G between two distinct vertices x and y is a sequence x, $z_1$, . . ., $z_k$, y of distinct vertices in G such that (i) x is adjacent to $z_1$; (ii) y is adjacent to $z_k$; (iii) if {$z_1$, . . ., $z_k$} ≠ ϕ then $z_i$ is adjacent to $z_{i+1}$ for i = 1 through k − 1 and (iv) if {$z_1$, . . ., $z_k$} = ϕ then x is adjacent to y. Obviously, if x and y are adjacent, then the sequence x, y is a path between them, with {$z_1$, . . ., $z_k$} = ϕ for this path. A vertex x is *connected to* a vertex y if there is a path between x and y. G is a *connected graph* if x is connected to y whenever x and y are distinct vertices of G.

### 2.4: Independent set and vertex cover

If I ⊂ V then I is an *independent set* of G if no two vertices in I are adjacent. If S ⊂ V then S *covers* an edge {x, y} if S ∩ {x, y} ≠ ϕ - that is, if S contains at least one end point of {x, y}. S is a *vertex cover* (of G) if S covers every edge of G. Additionally, if no proper subset of S is a vertex cover, then S is a *minimal vertex cover*. S is a *minimum vertex cover* if (i) S is a vertex cover and (ii) $|S| \neq |T|$ whenever T is a vertex cover of G. If S is a minimum vertex cover, then the positive integer $|S|$ is the *vertex cover number* (denoted by β(G)) of G. The following are immediate:
(i) S is a vertex cover if and only if V − S is independent;
(ii) S is a minimal vertex cover if and only if V − S is a maximal independent set; and
(iii) S is a minimum vertex cover if and only if V − S is a maximum independent set.

### 2.5: Recalling two essential theorems

**Theorem 2.1 - First theorem of graph theory** [12]: $\sum_{x \in G} dx = 2|E|$ in any graph G.

**Theorem 2.2 - Lower bounds for the vertex cover number** [9]: Let G be any graph.
(i) If M is a maximal matching of G, then $|M| \leq \beta(G)$; and
(ii) if $C_1, \ldots, C_k$ is any clique partition of G, then
$\}: \sum_{i=1 \text{ to } k} (|C_i - 1|) = n - k \leq \beta(G)$.

## 3. Further theory essential to the proposed algorithm

Throughout this section, G = (V, E) unless explicitly changed.

**Proposition 3.1.** A is an independent set of G if and only if dx(A) = 0 for all x ∈ A.

**Proof.** Follows from the definitions of independent set and dx(A) seen in subsection 2.2. ∎

**Proposition 3.2.** If I is an independent set of G then $\sum_{x \in I} dx \leq |E|$.

**Proof.** Let $F = \{X \in E : X \cap I \neq \phi\}$ and $W = \bigcup_{X \in F} X$. Clearly W − I is an independent set in the subgraph J = (W, F), and so J is a bipartite subgraph of G with the vertex set partition W = I ∪ (W − I). Further, $|F| \leq |E|$. Next, $\sum_{x \in I} dx = |F|$ in the bipartite graph J, and the conclusion follows. ∎

**Proposition 3.3.** If S is a vertex cover of G, then $\sum_{x \in S} dx \geq |E|$.

**Proof.** V − S is independent in G. By theorem 2.1, $\sum_{x \in S} dx + \sum_{x \in V - S} dx = 2|E|$. By proposition 3.2, $\sum_{x \in V - S} dx \leq |E|$, whence $\sum_{x \in S} dx \geq |E|$. ∎

**Proposition 3.4.** Let Δ = max {dy: y ∈ V} where dy denotes the degree (in G) of the vertex y. Then for any vertex cover S of G, $|S| \geq cf[|E|/\Delta]$ where for a real number m, cf [m] denotes the smallest integer k such that k ≥ m. (The function cf [m] is the *ceiling function* at m.)

**Proof.** Let S be a given vertex cover. Write $S = \{y_1, \ldots, y_p\}$. Let q and θ be, respectively, the integer part and the fractional part of $|E|/\Delta$, so that $|E|/\Delta = q + \theta$, with 0 ≤ θ < 1. By proposition 3.3, $\sum_{i=1 \text{ to } p} dy_i \geq |E|$.

If θ = 0 then $|E| = \Delta q$ and $cf[|E|/\Delta] = q$. Were $p (=|S|) < cf[|E|/\Delta] (= q)$, then $|E| \leq \sum_{i=1 \text{ to } p} dy_i \leq \Delta p < \Delta q = |E|$. Contradiction.
If θ > 0 then $|E| = \Delta q + \theta q$. Were $p (=|S|) < cf[|E|/\Delta] (= q + 1)$, then $|E| \leq \sum_{i=1 \text{ to } p} dy_i \leq \Delta p \leq \Delta q = |E| - \theta q < |E|$. Contradiction. ∎

**Corollary to 3.4.** $\beta(G) \geq cf[|E|/\Delta]$.

**Proposition 3.5.** Let A ⊂ V such that dx ≥ dy whenever x ∈ A and y ∈ V − A. If $\sum_{x \in A} dx < |E|$ and if B is a vertex cover of G then $|B| > |A|$.

**Proof.** By proposition 3.3, $\sum_{x \in B} dx \geq |E|$ and A is not a vertex cover. Assume $|B| \leq |A|$. Then neither $A \subset B$ nor $B \subset A$. If $A \cap B = \phi$ then $|E| > \sum_{x \in A} dx \geq \sum_{y \in B} dy \geq |E|$. Contradiction.

If $A \cap B \neq \phi$ then $|B - A| \leq |A - B|$ and so $\sum_{x \in B - A} dx \leq \sum_{x \in A - B} dx$. Then $|E| \leq \sum_{x \in B} dx = \sum_{x \in B - A} dx + \sum_{x \in B \cap A} dx \leq \sum_{x \in A - B} dx + \sum_{x \in B \cap A} dx = \sum_{x \in A} dx < |E|$. Contradiction. ∎

**Proposition 3.6.** Suppose S is a vertex cover of G. Then S is a minimal vertex cover if and only if for each $x \in S$ there is an edge X of G such that $S \cap X = \{x\}$.

**Proof.** ($\rightarrow$) Suppose S is a minimal vertex cover. Assume there is $y \in S$ for which the conclusion does not follow. Then $N(y) \subset S - \{y\}$. Let $D = S - \{y\}$. Let $\{a, b\} \in E$ be arbitrary. Then $a \in S$ or $b \in S$. If neither a nor b is y, then $\{a, b\} \cap D \neq \phi$. If either $a = y$ or $b = y$, then $\{a, b\} = \{y, z\}$ for some $z \in N(y)$, and in view of $N(y) \subset D$ it is immediate that $\{a, b\} \cap D \neq \phi$. Thus D is a vertex cover. Contradiction.
($\leftarrow$) Assume each $x \in S$ has the specified property. Let $y \in S$ be given. Then, by the hypothesis, there is a vertex $z \in V - S$ such that $\{y, z\} \in E$. Then the set $S - \{y\}$ does not cover the edge $\{y, z\}$. Consequently S is a minimal vertex cover. ∎

**Corollary to 3.6.** Let S be a vertex cover of G. Then S is a minimal vertex cover if and only if $dx(V - S) > 0$ for all $x \in S$.

**Proposition 3.7.** If $G = (V, E)$ is bipartite with the vertex set partitioning $V = A \cup B$ then $\min \{|A|, |B|\}$ is an upper bound on the vertex cover number of G.

**Proof.** Both A and B are vertex covers of G. ∎

## 4. The proposed Enter-Exit Greedy Algorithm (EEGA)

The crucial first step to fix a lower bound for the vertex cover number is proposition 3.3. Phase 1 of the EEGA is built on propositions 3.2 through 3.5. Phase 2 fixes an upper bound, and is built on the idea of dx(A) (subsection 2.2), propositions 3.1, 3.6 and 3.7.

### OUTLINE OF THE EEGA

**Input:** The vertex set V, the adjacency list and the adjacency matrix of a connected graph G.

### Phase 1: Lower bound for the vertex cover number of G

**Step 1.** (1a) $n = |V|$; go to (1b)
 (1b) Compute the degree of each $x \in V$ and order the vertices as:
 $x_1, \ldots, x_n$ where $dx_i \geq dx_{i+1}$ for $i = 1$ through n - 1; go to step 2
**Step 2.** (2a) Set $dx_i = d_i$ for $i = 1$ through n and go to (2b)
 (2b) Find the smallest positive integer L such that $d_1 + \ldots + d_L \geq |E|$; go to (2c)
 (2c) Return L and go to step 3

**Phase 2: Upper bound for the vertex cover number of $G$** \\

**Step 3:** $S = \{x_1, \ldots, x_L\}$ where $x_1$ through $x_L$ are the first L vertices in the ordering of V seen in (1b). Go to step 4

**Step 4:** If there is $x \in V - S$ such that the updated $d_x(V - S) > 0$
then choose the first such x (in the order seen in (1b)) and go to step 5,
else go to step 6

**Step 5:** Do updates: (i) $V - S \leftarrow V - S - \{x\}$ and (ii) $S \leftarrow S \cup \{x\}$
and return to step 4

**Step 6:** If there is $y \in S$ such that the updated $d_y(V - S) = 0$
then choose the last such y (in the order of (1b)) and go to step 7,
else go to step 8

**Step 7:**} Do updates: (i) $S \leftarrow S - \{y\}$ and (ii) $V - S \leftarrow (V - S) \cup \{y\}$
and return to step 6

**Step 8:** (8a) If $d_x(S) = 0$ for each $x \in S$
then $U = \min\{|S|, n - |S|\}$
else $U = |S|$ ; go to (8b)
(8b) If G is bipartite go to step 9, else go to step 10

**Step 9:** If $V = A \cup B$ is a partition of V
then $U \leftarrow \min\{U, |A|, |B|\}$ and go to step 10

**Step 10:** Return U

**Output:** $L \leq \beta(G) \leq U$.

---

In propositions 4.1 through 4.4, the input graph is $G = (V, E)$. \\

**Proposition 4.1.** The EEGA is feasible - that is, it terminates in finitely many computations.

**Proof.** Clearly $d_1 + \ldots + d_k \geq |E|$ for $k = n - 1$. So there is a smallest positive integer L such that $d_1 + \ldots + d_L \geq |E|$. Hence phase 1 terminates in finitely many computations.
There are only finitely many $x \in V - S$ with $d_x(V - S) > 0$. Transfering each such x to S from $V - S$ (see step 5 in phase 2) is feasible. Next, the exits from S to $V - S$ (see step 7 in phase 2) are all feasible because there are only finitely many $y \in S$ with $d_y(V - S) = 0$ and no such y is transferred back to S. Thus phase 2 of the EEGA also terminates in finitely many computations. ∎

**Proposition 4.2.** The positive integer L returned in phase 1 of the EEGA is a lower bound for $\beta(G)$.

**Proof.** Note that L is the smallest positive integer returned by (2c). Suppose S is a vertex cover with $|S| = q$. Then $d_1 + \ldots + d_q \geq |E|$ by proposition 3.3, and so the integer L returned by (2c) cannot be larger than q. Hence $L \leq |S|$ for any vertex cover S. ∎

**Proposition 4.3.** The EEGA converges to a desired output.

**Proof.** A desired output is a minimal vertex cover. When the algorithm enters step 6, $d_x(V - S) = 0$ for each $x \in V - S$. By proposition 3.1, $V - S$ is an independent set. So S is a vertex cover. Next, when the algorithm enters step 8, $d_x(V - S) > 0$ for each $x \in S$ and the control does not return to any preceding step. By corollary to 3.6, S is a minimal vertex cover. Hence the EEGA converges to a minimal vertex cover. ∎

**Proposition 4.4.** Phase 2 returns a positive integer that is no less than L.

**Proof.** When phase 2 begins, $|S| = L$ where L is (the lower bound) returned by phase 1. So, if $B \subset V$ and $|B| < L$ then clearly $\sum_{x \in B} d_x < |E|$ and so B cannot be a vertex cover. Steps 5 through 8 ensure that either S or $V - S$ is a minimal vertex cover (by dint of proposition 3.6), and that this minimal vertex cover is of size U. Hence $U \geq L$. ∎

## 5. The worst-case time complexity of the EEGA

The worst-case time complexiy of the EEGA is discussed using the asymptotic growth rate functions big-oh (O) [8, 17, 22]. By default, the term time complexity will mean the worst-case one, throughout this section. The following property [17] is important to the discussion in this section. This property will be referred to as k-sum property, and will be invoked in the last subsection of this section.

*k-sum property:* Let k be a fixed positive integer (that is, k does not depend on any input size in the algorithm under discussion). If $f_i$ (i = 1, . . ., k) and h are functions such that $f_i = O(h)$ for all i = 1, . . ., k, then $f_1 + \ldots + f_k = O(h)$.

In section 4, the algorithm is presented in pseudo-code style, in steps 1 through 10. In this section we show the EEGA has a polynomial running time of $O(n^3)$ (where n is the order of the input graph). The pseudo-code has been divided into nine parts for analysis, and each part is presented in a subsection. From the pseudo-code of the EEGA, it is clear that these parts run in sequence. The input graph G = (V, E) is represented by its vertex set and its adjacency list.

Throughout this section, by the phrase "(*) is bounded by time $O(n^d)$," we will mean that there are absolute constants $c > 0$ and $d > 0$ so that on every input graph of order n, the running time of the process in the place of (*) is bounded by $cn^d$ primitive computational steps ([17], chapter 2). The following are the primitive computational steps in the EEGA:
(p-c 1) Assigning a value to a variable
(p-c 2) Placing a new element at the end of a list of elements
(p-c 3) Reading an element from a list
(p-c 4) Any of the four fundamental operation on real numbers
Further, the term "instance" in this section will mean an input graph.

### 5.1: Part 1: Step 1 of the pseudocode - time complexity $T_1$

The computation of the degrees of the vertices of G can be done in $O(n|E|)$ time. The time complexity for sorting the degrees in non-ascending order is $O(n^2)$. In each instance, step 1 is run once. Hence $T_1 = O(n^2)$.

### 5.2: Part 2: Step 2 - time complexity $T_2$

(2a) consists of n assignments, each of constant time (that is, O(1)). In (2b), the process of computing L involves the cumulative sums $d_1 + \ldots + d_r$ ($1 \leq r \leq n - 1$). This clearly does not exceed n − 1 additions. (2c) obviously takes only constant time. In each instance, step 2 is run once. Hence $T_2 = O(n)$.

### 5.3: Part 3: Step 3 - time complexity $T_3$

In step 3, the set $S = \{x_1, \ldots, x_L\}$ is obtained by beginning with the singleton $\{x_1\}$ and queueing the remaining L − 1 vertices one by one after $x_1$. Queueing each vertex behind its predecessors takes constant time [17]. In each instance, step 3 is run once. Hence $T_3 = O(n)$.

### 5.4: Part 4: Step 4 - time complexity $T_4$

Finding the first $x \in V - S$ such that $dx(V - S) > 0$ is bounded by $|V - S|$ readings from the adjacency list, and so is bounded by time $O(n^2)$. In each instance, step 4 is run once every time the control returns to the step from step 5. But $|V - S| \leq n - 1$ when step 5 begins, and each update $V - S \leftarrow V - S - \{x\}$ in this step only reduces $|V - S|$. So step 5 returns the control to step 4 at most n − 2 times. Hence in each instance, step 4 is run not more than n times. Hence $T_4 = O(n^3)$.

### 5.5: Part 5: Step 5 - time complexity $T_5$

Each update $V - S \leftarrow V - S - \{x\}$ takes O(n) time, and there are at most n − 2 such updates. Next, each update $S \leftarrow S \cup \{x\}$ is done in constant time because it is a queueing operation, and there are at most n − 2 such updates. Also, there are only finitely many such $x \in V - S$ that would be deleted from V − S and would be added to S, which means the control goes to step 6 in finitely many computations and never returns to step 5. So, in each instance, step 5 is run not more than n times. Hence $T_5 = O(n^2)$.

### 5.6: Part 6: Step 6 - time complexity $T_6$

Finding the last $y \in S$ such that $dy(V - S) = 0$ requires reading at most $|S|$ values, and so is bounded by $O(n^2)$. In each instance, step 6 is run once every time the control returns to the step from step 7. But then $|S| \leq n - 1$ when step 7 begins, and each update $S \leftarrow S - \{y\}$ in this step only reduces $|S|$. So step 7 returns the control to step 6 at most n − 2 times. Hence in each instance, step 6 is run not more than n times, from which $T_6 = O(n^3)$.

### 5.7: Part 7: Step 7 - time complexity $T_7$

The operations in step 7 are similar to those in step 5. Hence $T_7 = O(n^2)$.

### 5.8: Part 8: Steps 8 - time complexity $T_8$

Checking if $dx(S) = 0$ requires $O(n)$ time for each $x \in S$. Assignation of the right value to the variable U is done in constant time. As for (8b), it is done using breadth-first search, in $O(n^2)$ time. Hence $T_8 = O(n^2)$.

### 5.9: Part 9: Steps 9 and 10 - time complexity $T_9$

Obviously $|A| \leq n - 1$ and $|B| \leq n - 1$. Step 9 determines $|A|$ and $|B|$, and this is bounded by $O(n)$ time. This is followed by determining the least of three positive integers and assigning it to the variable U, which is done in constant time. Hence $T_9 = O(n)$.

### 5.10: Time complexity T of the EEGA

The EEGA has no additional parts other than the nine shown in the preceding nine parts (subsections 5.1 through 5.9). Also, from the preceding subsections, it follows that $T_i = O(n^3)$ for each $i = 1, \ldots, 9$ because if $f(n) = O(n^k)$ then $f(n) = O(n^r)$ for all $r > k$ ([17], chapter 2). Here we invoke the k-sum property placed in the beginning of this section, with $k = 9$ (which number does not increase in the algorithm). Hence the worst-case time complexity of the EEGA is $T = \sum_{i=1 \text{ to } 9} T_i = O(n^3)$.

## 6. Discussion

The following observation is at the base of the EEGA: if $G = (V, E)$ is a graph and $S \subset V$ with $\sum_{x \in S} dx < |E|$ then S cannot be a vertex cover of G.

### 6.1: Comparing the EEGA with a classical maximal matching algorithm

A classical algorithm (see (i) of theorem 2.2 in section 2 for the basic fact) finds a maximal matching M and fixes $|M|$ as a lower bound for $\beta(G)$. The vertex cover output by this matching algorithm is of size $2|M|$ since all the vertices of M are returned. This output is not necessarily a minimal vertex cover.
In contrast, the EEGA, though it might return a lower bound less than $|M|$, always returns a minimal vertex cover. This feature ensures that the upper bound returned by the EEGA is at least as good as that returned by the classical maximal matching algorithm for any graph. Also, in many an instance this feature ensures that the EEGA returns the precise $\beta(G)$.

### 6.2: Comparing the EEGA with a classical clique partition algorithm

A classical clique partition algorithm (see (ii) of theroem 2.2 in section 2 for the basic fact) partitions the vertex set into cliques to fix a lower bound for $\beta(G)$. The vertex cover output by this algorithm is the union of all the non-trivial cliques in this partition [9]. This, like in the classic maximal matching algorithm mentioned above, is not necessarily minimal. (Sometimes the output is the whole of V.)
In contrast, the EEGA ensures that a minimal vertex cover is always returned. So the upper bound returned by the EEGA is at least as good as that returned by the classical clique partition algorithm for any graph.

### 6.3: Comparing the EEGA with greedy algorithms

The EEGA has provisions (in phase 2) for deletion of vertices from existing vertex cover without affecting the covering property. This is what ensures that the upper bound for β(G) arises from a minimal vertex cover. In this aspect the EEGA is an improvement upon the existing greedy algorithms.

### 6.4: Other merits of the EEGA

(i) It runs in $O(n^3)$ time.
(ii) For any tree, it always outputs a minimal vertex cover not containing any leaf. In particular, it returns the optimal solution for any path graph.
(iii) The initial ordering of the vertices in non-ascending order of their degrees helps in making a good choice of a vertex cover `candidate' set to begin phase 2. It is from this set that a minimal vertex cover is constructed.

### 6.5: Limitations of the EEGA

(i) Any tie in step 4 (phase 2) is broken by choosing the first vertex on the tied list (in the order seen in (1b)). This may result in the algorithm missing a better (that is, numerically lesser) upper bound.
(ii) The algorithm fixes bounds but does not give any indication as to which of them is closer to the optimal solution.
(iii) The EEGA does not guarantee the optimal solution except when it returns equal L and U.

### 6.6: An illustration

Let G = (V, E) be the graph with vertex set V = {a, b, c, d, e, f, g, h, i, j, x, y, z}. Also, the adjacency list of G is:
(i) N(a) = {b, c}; (ii) N(b) = {a, e, f, i}; (iii) N(c) = {a, d, f, g};
(iv) N(d) = {c, j}; (v) N(e) = {b, h}; (vi) N(f) = {b, c, i, j};
(vii) N(g) = {c, i, j}; (viii) N(h) = {e, y};
(ix) N(i) = {b, f, g, h, x}; (x) N(j) = {d, f, g, x, z};
(xi) N(x) = {i, j, y, z}; (xii) N(y) = {h, x, z} and (xiii) N(z) = {j, y, x}.
Here n = |V| = 13 and |E| = 22.

### EEGA: Phase 1
(1) Let $A(\Delta_i)$ = {x ∈ V : dx = i}. Then V = $A(\Delta_5)$ ∪ $A(\Delta_4)$ ∪ $A(\Delta_3)$ ∪ $A(\Delta_2)$,
  where $A(\Delta_k)$ ∩ $A(\Delta_p)$ = φ for k ≠ p.
(2) Here $A(\Delta_5)$ = {i, j}; $A(\Delta_4)$ = {b, c, f, x}; $A(\Delta_3)$ = {g, h, y, z} and
  $A(\Delta_2)$ = {a, d, e}.
 (3) Ordering of the 13 vertices of G: i, j, b, c, f, x, g, h, y, z, a, d and e
    (by (1b) of the EEGA).
 (4) Their degrees are: 5, 5, 4, 4, 4, 4, 3, 3, 3, 3, 2, 2 and 2, respectively.
    These degrees are designated $d_1$ through $d_{13}$, respectively.
    Also, Δ = 5 and cf[ |E| / Δ] = 5.

(5) From (4), the smallest positive integer L such that $\sum_{i=1 \text{ to } L} d_i \geq |E|$ is L = 5.
   (That is, 5 + 5 + 4 + 4 + 4)
(6) Hence $\beta(G) \geq 5$.

### EEGA: Phase 2
(7) Initial vertex cover candidate is S = {i, j, b, c, f}, the set of the first L (= 5)
   vertices of the ordered list in (3) of phase 1. Then V − S = {x, g, h, y, z, a, d, e}.
(8) Entry of vertices into S (from V − S)
   (8a) u = x is the first vertex in V − S for which du (V − S) > 0.
   Updates: S ← S ∪ {x} and V − S arrow V − S - {x}
   Now S = {i, j, b, c, f, x} and V − S = {g, h, y, z, a, d, e}
   (8b) u = h is the first vertex in V − S for which du(V − S) > 0.
   Updates: S ← S ∪ {h} and V − S ← V − S - {h}
   Now S = {i, j, b, c, f, x, h} and V − S = {g, y, z, a, d, e}
   (8c) u = y is the first vertex in V − S for which du(V − S) > 0.
   Updates: S ← S ∪ {y} and V − S ← V − S - {y}
   Now S = {i, j, b, c, f, x, h, y} and V − S = {g, z, a, d, e}
   Here V − S is independent since du(V − S) = 0 for all u ∈ V − S.
   So no more entry of vertices into S.
(9) Exit of vertices from S (to V − S)
   (9a) u = f is the last vertex in S for which du(V − S) = 0.
   Updates: S ← S - {f} and V − S ← (V − S) ∪ {f}
   Now: S = {i, j, b, c, x, h, y} and V − S = {f, g, z, a, d, e}
   Here S is a minimal vertex cover since du(V − S) > 0 for all u ∈ S.
   So further vertex exit from S is not allowed.
(10) Since S = {i, j, b, c, x, h, y} is not an independent set and G is not bipartite,
   return U = |S| = 7.
(11) **Output: $5 \leq \beta(G) \leq 7$**.

## 7. Concluding remarks

We have presented an efficient algorithm (EEGA) that fixes a lower bound L and an upper bound U for the vertex cover number of a given graph. We wish to emphasize that for many a graph the EEGA gives the optimal solution by returning U = L.
We also emphasize that the EEGA is guaranteed to return a minimal vertex cover. This is not so with any of the earlier greedy algorithms.
The vertices are input in non-ascending order of their degrees, beginning with one having the largest degree in the given graph G. The results of section 3 show this is a good ordering for phase 1 of the EEGA. But we do not, at this point, know if this order would work well in phase 2 for all graphs.
Our experiments with benchmark graphs of various sizes (up to 2000 vertices) led us to optimal solutions for many graphs (meaning, L = U). But for many other graphs the EEGA could only output U larger than L, thereby only building intervals around their optimum solutions.
Since the EEGA fixes bounds on β(G) and is not (yet) an approximation algorithm for β(G), it can be compared only with algorithms that fix bounds, as has been done in subsections 6.1 through 6.3. So we have not reported our experiments in detail.
Further, for many graphs the EEGA returned L and U that are `close' to each other, meaning U − L is small compared to L. The authors are currently working more on the EEGA, exploring for possibilities of an additional polynomial-time phase that could conclusively report on the possibility of vertex covers of size smaller than U.

## Acknowledgements

This research was supported fully by Niels Abel Foundation, Palakkad, Kerala State, India.

## References

[1] Akiba, T., Iwata, Y.: Branch-and-reduce exponential/fpt algorithms in practice: A case study of vertex cover. CoRR, abs/1411.2680 (2014).

[2] Angel, E., Campigotto, R., Laforest, C.: Algorithms for the vertex cover problem on large graphs. IBISC – Universit´e d'Evry-Val d'Essonne Research report no. 2010-01 (2010).

[3] Balakrishnan, R., Ranganathan, K.: A textbook of Graph Theory. 2nd edn. Springer, New York (2012).

[4] Berge, C.: Graphs and Hypergraphs. North-Holland, Amsterdam (1973).

[5] Bourgeois, N., Escoffier, B., Paschos, V.T., van Rooij, J.M.M.: Fast Algorithms for max independent set. Algorithmica 62}(1), 382-415 (2012).

[6] Chen, T., Kanj, I.A., Xia, G.: Improved upper bounds for vertex cover. Theoretical Computer Science 411}, 3736-3756 (2010).

[7] Cook, S.: The P versus NP problem. Available online: http://www.claymath.org/millennium/P vs NP/pvsnp.pdf

[8] Cormen, T.H., Leiserson, C.E., Rivest, R.L., Stein, C.: Introduction to Algorithms. 3rd edn. MIT Press, Cambridge (2009).

[9] Delbot, F., Laforest, C., Phan, R.: New Approximation Algorithms for the vertex cover problem and variants. Research Report LIMOS/RR-13-02 (2013).

[10] Deo, N.: Graph theory with applications to Engineering and Computer Science. Prentice Hall of India, New Delhi (1984).

[11] Garey, M.R., Johnson, D.S.: Computers and Intractability: A guide to the Theory of NP-completeness. Freeman, New York (1979).

[12] Harris, J. M., Hirst, J.L., Mossinghoff, M.J.: Combinatorics and Graph Theory. 2nd edn. Springer, New York (2008).

[13] Islam, A-U., Kalita, B.: Application of Minimum Vertex Cover for Keyword – based Text Summarization Process. International Journal of Computational Intelligence Research 13}(1), 113-125 (2017).

[14] Kanj, I.: Vertex cover - Exact and Approximation Algorithms and Appplications. Ph. D. Thesis, Texas A \& M University, United States (2001).

[15] Karakostas, G.: A better approximation ratio for the vertex cover problem. In: ICALP 2005, Lisbon, Portugal, LNCS, vol. 3580, pp. 1043-1050. Springer-Verlag, Heidelberg (2005).

[16] Karp, R.: Reducibility among Combinatorial Problems. In: Miller, R., Thatcher, J. (eds.) Complexity of Computer Computations, pp. 85-103. Plenum Press, New York (1972).

[17] Kleinberg, J., Tardos, E.: Algorithm Design. Pearson Addison Wesley, New York (2006).

[18] Papadimitriou, C., Yannakakis, M.: Optimization, Approximation, and Complexity Classes. In: STOC 88: Proceedings of the 20th annual ACM symposium on Theory of Computing, pp. 229-234. ACM Press, New York (1988).

[19] Savage, C.: Depth-first search and the vertex cover problem. Information Processing Letters 14}(5), 233-235 (1982).

[20] Stoll, R.R.: Set Theory and Logic. Dover Publications, New York (1963).

[21] Urhaug, T.S.: A Survey of Linear-Programming Guided Branching Parameterized Algorithms for Vertex Cover with Experimental Results. Master Thesis in Computer Science, Department of Informatics, University of Bergen, Norway (2015).

[22] Wilf, H.S.: Algorithms and Complexity. Internet Edition, Summer (1994), available on http://www/cis.upenn.edu/wilf